\PassOptionsToPackage{unicode}{hyperref}
\PassOptionsToPackage{hyphens}{url}
\PassOptionsToPackage{dvipsnames,svgnames,x11names}{xcolor}
\documentclass[
  12pt]{article}

\usepackage{amsmath,amssymb,mathtools}
\usepackage{bm}
\usepackage{iftex}
\ifPDFTeX
  \usepackage[T1]{fontenc}
  \usepackage[utf8]{inputenc}
  \usepackage{textcomp} %
\else %
  \usepackage{unicode-math}
  \defaultfontfeatures{Scale=MatchLowercase}
  \defaultfontfeatures[\rmfamily]{Ligatures=TeX,Scale=1}
\fi
\usepackage{lmodern}
\usepackage{algorithm}
\usepackage{algpseudocode}
\ifPDFTeX\else
\fi
\IfFileExists{upquote.sty}{\usepackage{upquote}}{}
\IfFileExists{microtype.sty}{%
  \usepackage[]{microtype}
  \UseMicrotypeSet[protrusion]{basicmath} %
}{}
\makeatletter
\@ifundefined{KOMAClassName}{%
  \setlength{\parindent}{2em}
  \setlength{\parskip}{0pt}
}{%
  \KOMAoptions{parskip=false}
  \setlength{\parindent}{2em}}
\makeatother
\usepackage{indentfirst}
\usepackage{xcolor}
\makeatletter
\ifx\paragraph\undefined\else
  \renewcommand{\paragraph}{\@ifstar\xxxParagraphInline\xxxParagraphInline}
  \newcommand{\xxxParagraphInline}[1]{\par\noindent\emph{#1}\ }
\fi
\ifx\subparagraph\undefined\else
  \renewcommand{\subparagraph}{\@ifstar\xxxSubParagraphInline\xxxSubParagraphInline}
  \newcommand{\xxxSubParagraphInline}[1]{\par\noindent\emph{#1}\ }
\fi
\makeatother

\usepackage{longtable,booktabs,array}
\usepackage{calc} %
\usepackage{etoolbox}
\makeatletter
\patchcmd\longtable{\par}{\if@noskipsec\mbox{}\fi\par}{}{}
\makeatother
\IfFileExists{footnotehyper.sty}{\usepackage{footnotehyper}}{\usepackage{footnote}}
\makesavenoteenv{longtable}
\usepackage{graphicx}
\usepackage{placeins}
\makeatletter
\def\maxwidth{\ifdim\Gin@nat@width>\linewidth\linewidth\else\Gin@nat@width\fi}
\def\maxheight{\ifdim\Gin@nat@height>\textheight\textheight\else\Gin@nat@height\fi}
\makeatother
\setkeys{Gin}{width=\maxwidth,height=\maxheight,keepaspectratio}
\makeatletter
\def\fps@figure{htbp}
\makeatother
\makeatletter
\@ifpackageloaded{caption}{}{\usepackage{caption}}
\AtBeginDocument{%
\ifdefined\contentsname
  \renewcommand*\contentsname{Table of contents}
\else
  \newcommand\contentsname{Table of contents}
\fi
\ifdefined\listfigurename
  \renewcommand*\listfigurename{List of Figures}
\else
  \newcommand\listfigurename{List of Figures}
\fi
\ifdefined\listtablename
  \renewcommand*\listtablename{List of Tables}
\else
  \newcommand\listtablename{List of Tables}
\fi
\ifdefined\figurename
  \renewcommand*\figurename{Figure}
\else
  \newcommand\figurename{Figure}
\fi
\ifdefined\tablename
  \renewcommand*\tablename{Table}
\else
  \newcommand\tablename{Table}
\fi
}
\@ifpackageloaded{float}{}{\usepackage{float}}
\floatstyle{ruled}
\@ifundefined{c@chapter}{\newfloat{codelisting}{h}{lop}}{\newfloat{codelisting}{h}{lop}[chapter]}
\floatname{codelisting}{Listing}

\makeatother
\makeatletter
\@ifpackageloaded{caption}{}{\usepackage{caption}}
\@ifpackageloaded{subcaption}{}{\usepackage{subcaption}}
\makeatother

\ifLuaTeX
  \usepackage{selnolig}  %
\fi
\usepackage[authoryear,round]{natbib}
\usepackage{bookmark}
\usepackage{etoc}

\IfFileExists{xurl.sty}{\usepackage{xurl}}{} %
\hypersetup{
  pdftitle={Human-Anchored Inference for Ranking New Models with Large Language Model Judges},
  pdfauthor={Anonymous},
  pdfkeywords={Bradley-Terry-Luce model; covariate shift; data fusion; human-reference ranking; LLM judge; pairwise comparison; semiparametric inference},
  colorlinks=true,
  linkcolor={blue},
  filecolor={Maroon},
  citecolor={Blue},
  urlcolor={Blue},
  pdfcreator={LaTeX via pandoc}}

\usepackage{amsthm}

\usepackage{booktabs,array,enumitem,microtype}
\newcommand{\anon}{1}
\if1\anon
  \hypersetup{pdfauthor={Xin Zhou, Sinian Zhang, Zhanyan Yang, Mingyuan Xu, Molei Liu, Doudou Zhou}}
\fi

\newtheorem{theorem}{Theorem}[section]
\newtheorem{corollary}{Corollary}[section]

\newtheorem{assumption}{Assumption}[section]
\theoremstyle{definition}

\newtheoremstyle{italicremark}
  {3pt}{3pt}
  {\itshape}
  {}
  {\itshape}
  {.}
  {.5em}
  {}
\theoremstyle{italicremark}

\newcommand{\E}{\mathbb{E}}
\newcommand{\Prob}{\mathbb{P}}

\newcommand{\Hist}{\mathrm{Hist}}

\newcommand{\Hum}{\mathrm{H}}

\newcommand{\cf}{\mathrm{cf}}
\newcommand{\argmin}{\mathop{\mathrm{argmin}}}

\newcommand{\sigmoid}{\sigma}

\newcommand{\vect}[1]{\bm{#1}}
\newcommand{\mat}[1]{\bm{#1}}

\begin{document}
\etocdepthtag.toc{main}

\def\spacingset#1{\renewcommand{\baselinestretch}%
{#1}\small\normalsize} \spacingset{1}

\if1\anon
{
  \title{\bf Human-Anchored Inference for Ranking New Models with Large Language Model Judges}
  \author{%
    Xin Zhou\textsuperscript{1,*}, Sinian Zhang\textsuperscript{2,*},
    Zhanyan Yang\textsuperscript{1},\\
    Mingyuan Xu\textsuperscript{3}, Molei Liu\textsuperscript{4,\textdagger},
    Doudou Zhou\textsuperscript{3,\textdagger}\\[1em]
    \begin{minipage}{0.95\textwidth}
      \centering\small
      \textsuperscript{1}School of Management, University of Science and Technology of China\\[0.5em]
      \textsuperscript{2}Division of Biostatistics and Health Data Science, University of Minnesota\\[0.5em]
      \textsuperscript{3}Department of Statistics and Data Science, National University of Singapore\\[0.5em]
      \textsuperscript{4}Department of Biostatistics, Peking University Health Science Center, Peking University\\[0.75em]
      \textsuperscript{*}Xin Zhou and Sinian Zhang contributed equally to this work.\\[0.3em]
      \textsuperscript{\textdagger}Co-corresponding authors:\\[0.2em]
      Molei Liu (\href{mailto:moleiliu@bjmu.edu.cn}{moleiliu@bjmu.edu.cn});
      Doudou Zhou (\href{mailto:ddzhou@nus.edu.sg}{ddzhou@nus.edu.sg}).
    \end{minipage}}
  \date{}
  \maketitle
} \fi

\if0\anon
{
  \bigskip
  \bigskip
  \bigskip
  \begin{center}
    {\LARGE\bf Human-Anchored Inference for Ranking New Models with Large Language Model Judges}
\end{center}
  \medskip
} \fi

\bigskip

\begin{abstract}
Human pairwise comparisons provide a reference for evaluating large language models (LLMs), but collecting sufficient judgments for each new release is costly and time-consuming. LLM judges offer a scalable alternative, although their comparisons may differ systematically from human preferences and across judges. We study the ranking of a new model that has received LLM-judge comparisons but no human comparisons. We propose ANCHOR (ANchored Comparisons for Human-reference inference with Orthogonal Riesz correction), which uses historical human and LLM comparisons to learn judge-specific sensitivities to human score differences and feature-dependent judge biases. These estimates are then used to infer the new model's human-reference score based on its LLM-judge comparisons only. The framework allows the feature distribution to change between historical and new-model comparisons. For inference, we construct a Neyman-orthogonal estimator through a joint Riesz correction that removes the first-order bias arising from the estimation errors of the historical human scores, judge sensitivities, and bias functions. We establish identification, convergence rates, and asymptotic normality. We also show that ANCHOR attains the semiparametric efficiency bound. Simulations demonstrate gains in score estimation and ranking accuracy. On Chatbot Arena, ANCHOR achieves the lowest score RMSE and insertion MAE among competing methods, with narrower score intervals on average.
\end{abstract}

\noindent%
\vfill

\newpage
\spacingset{1.8} %

\section{Introduction}

Pairwise rankings are often constructed from comparisons with heterogeneous evaluators or data sources \citep{fan2025heterogeneous}. In many applications, one evaluator population defines the target ranking but is costly to query, while alternative evaluators provide more abundant comparisons that may differ systematically from the reference preferences. A particular challenge arises when a new item has been evaluated only by these alternative sources and must nevertheless be ranked according to the reference population. Large language model (LLM) evaluation provides an instance of this problem. Frequent releases across major model families
create a continuing need to compare new models with existing ones. Human pairwise preferences provide a natural and interpretable reference for such comparisons, but collecting sufficient human judgments for every new release is costly and time-consuming \citep{chiang2024chatbot,zheng2023judging}. By contrast, LLM judges can generate comparisons much more quickly and at substantially lower cost \citep{li2025generation}. Their comparisons cannot be treated directly as substitutes for human preferences because different LLM judges may exhibit systematic biases and differ substantially in how they evaluate the same responses \citep{chen2024judgementbias}. The task is therefore to use historical human and LLM-judge comparisons to infer a newly released model's position in the human-reference ranking from its LLM-judge comparisons, before human comparisons involving that model become available.

A widely used model for ranking from pairwise comparisons is the Bradley--Terry--Luce (BTL) model \citep{bradley1952rank,luce1959individual}, which assigns each item a latent score and expresses comparison probabilities through score differences. A substantial literature has developed methods for score estimation and statistical inference under this model \citep{negahban2012iterative,chen2015spectral,chen2019spectral,simons1999asymptotics,liu2022lagrangian,gao2021uncertainty}. More recent work incorporates covariates, heterogeneous or contextual preferences, and distributional shift \citep{fan2024care,fan2025spectral,zhang2026fisher}. These methods primarily study ranking under a given comparison mechanism.

The growing use of LLMs as judges has also motivated methods that account for evaluator-specific behavior \citep{xu2026judgeaware,yu2026heterogeneous} and allow rankings to vary across prompts or tasks \citep{frick2025prompt,avelar2026prompt,li2026lowrankrank}. A related literature studies how LLM-judge comparisons can be combined with human comparisons. For example, \citet{chatzi2024ppr} use prediction-powered ranking (PPR) to construct rank sets, but their labeled sample requires every ranked model to appear in human comparisons. \citet{polo2025bridge} model systematic discrepancies between human and LLM judgments through observed response features and derive inference for bias-adjusted evaluation quantities, including predictions for new observations. \citet{guerdan2026doubly} combine human and LLM persona ratings for doubly robust evaluation under sampling bias, while \citet{chen2026efficient} derive efficient inference for human win rates. JudgeArena approximates LMArena Elo by fitting a common BTL model to historical human and target-model judge comparisons \citep{lushtaku2026judgearena}. The remaining problem is to infer a new model's score on the historical human BTL scale from heterogeneous LLM-judge comparisons. The feature distribution of these new-model comparisons may differ from that of historical LLM-judge comparisons.

This ranking problem creates two statistical difficulties. First, the new model has no human comparisons, and LLM judges may disagree systematically with humans. We therefore need to establish when historical human and LLM-judge comparisons can identify the new model's human-reference score, even when the feature distribution changes between historical and new-model comparisons. Second, the historical human scores and the relationship between human and LLM-judge comparisons must themselves be estimated. Their estimation errors propagate to the new-model score, so valid inference must account for both uncertainty in the historical estimates and sampling variation in the new-model comparisons.

To address these difficulties, we introduce a model linking human and LLM-judge comparisons through human-reference scores. Human comparisons follow a BTL model. For each LLM judge, the comparison logit combines a human score difference weighted by a judge-specific sensitivity with a bias determined by prompt and response features. Each judge's comparison mechanism is shared across historical and new-model comparisons, while their feature distributions may differ. Under this model, we establish conditions for identifying the new model's human-reference score from the observed comparisons.

Based on this model, we develop ANCHOR (ANchored Comparisons for Human-reference inference with Orthogonal Riesz correction). ANCHOR first estimates the historical scores from human comparisons, then fits each judge model conditional on these estimates. The fitted judge models provide an initial estimate of the new model's human-reference score from its LLM-judge comparisons. For inference, a joint Riesz correction combines weighted residuals from all three samples to yield a Neyman-orthogonal estimator. We establish convergence rates for the initial estimates and asymptotic linearity and normality of the corrected estimator. We further derive the semiparametric efficiency bound and show that ANCHOR attains it. {In specific, we identify and quantify an information bottleneck: increasing LLM-judge comparisons cannot remove uncertainty in the human score scale.} {For scores and contrasts affected by this scale, we show that ANCHOR's variance approaches a lower bound determined solely by historical human comparisons when the LLM-judge comparisons have a much larger sample size.} Simulations and empirical experiments on Chatbot Arena and MT-Bench assess finite-sample performance under judge heterogeneity and distribution shift.

The remaining paper is organized as follows. Section~\ref{sec:model} defines the data, comparison models, and identification conditions. Section~\ref{sec:method} introduces the two-step estimator and orthogonal correction, with theoretical results in Section~\ref{sec:theory-v2}. Sections~\ref{sec:simulation} and~\ref{sec:realdata} report simulation and empirical studies, and Section~\ref{sec:discussion-v2} concludes with a discussion. Proofs are provided in Supplementary Section~E, and auxiliary results are collected in Supplementary Section~F.

\section{Problem formulation}\label{sec:model}

\subsection{Observed comparison data}

Let \(\mathcal I=\{1,\ldots,K\}\) index the \(K\) historical models. We observe mutually independent samples of pairwise comparisons among these models from humans and from \(M\) LLM judges indexed by \(m\in\{1,\ldots,M\}\). The historical human and LLM-judge comparisons are
\begin{equation*}
\begin{aligned}
O_r^{\Hum}&:=(I_r^{\Hum},J_r^{\Hum},Y_r^{\Hum})\overset{\mathrm{i.i.d.}}{\sim}P_H,\\
O_r^{C_m}&:=(I_r^m,J_r^m,Y_r^{C_m},\vect{X}_{I_r^m,r}^{m},\vect{X}_{J_r^m,r}^{m})\overset{\mathrm{i.i.d.}}{\sim}P_{C_m}.
\end{aligned}
\end{equation*}
The indices \(I\) and \(J\) denote the models displayed on the left and right, and \(Y=1\) (\(Y=0\)) indicates a preference for the left (right) response. For an LLM-judge comparison \(r\), \(\vect X_{i,r}^{m}\in\mathbb R^{d_x}\), with \(i\in\{I_r^m,J_r^m\}\), contains features of model \(i\)'s displayed response and its prompt context, observed before judge \(m\)'s verdict. These features, such as response length, formatting indicators, and prompt--response representations, are used below to model systematic judge-specific deviations not explained by human score differences. After coordinatewise rescaling, we take \(\mathcal X=[-1,1]^{d_x}\) as the common feature space.

Let \(\star\notin\mathcal I\) denote a new model for which no human comparison is observed. Before observing the comparison outcomes, fix a nonempty set \(\mathcal M_A\subseteq\{1,\ldots,M\}\) of LLM judges used to evaluate it. The new model is compared with historical models only by these judges, producing a sample of \(n_\star\) comparisons independent of the historical samples
\begin{equation*}
O_s^A:=(\vect{X}_{\star,s}^{\star},\vect{X}_{J_s,s}^{\star},J_s,G_s,Y_s^A)
\overset{\mathrm{i.i.d.}}{\sim}P_A,
\qquad
\mathcal O_A:=\{O_s^A\}_{s=1}^{n_\star}.
\end{equation*}
In comparison \(s\), \(J_s\in\mathcal I\) is the historical model on the right, \(G_s\in\mathcal M_A\) is the judge index, and \(Y_s^A=1\) means that the judge prefers the new model's response. The superscript \(\star\) distinguishes features in new-model comparisons from those in historical comparisons. Each selected judge has positive sampling probability, so \(\Prob(G_s=m)>0\) for every \(m\in\mathcal M_A\).
For~judge \(m\in\mathcal M_A\), define
\(\mathcal O_{A_m}:=\{O_s^A:G_s=m\}\),
\(P_{A_m}:=P_A(\cdot\mid G=m)\), and
\(n_{A_m}:=\sum_{s=1}^{n_\star}\mathbf1\{G_s=m\}\).
Thus \(n_\star=\sum_{m\in\mathcal M_A}n_{A_m}\).
We index the human, historical judge, and new-model judge samples by
\(\mathcal S=\{H,C_1,\ldots,C_M\}\cup\{A_m:m\in\mathcal M_A\}\).
For \(q\in\mathcal S\), \(\mathcal O_q=\{O_r^q\}_{r=1}^{n_q}\) is its sample
of \(n_q\) comparisons from \(P_q\), and \(\E_q\) denotes expectation under
\(P_q\). Within each source, \(r\) indexes comparisons and \(Y_r^q\) is the
outcome; a subscript \(r\) on a function denotes its evaluation at that comparison.
We condition on the judge allocation when working with the samples
\(\mathcal O_{A_m}\). All sample sizes are indexed by a common asymptotic index \(n\), with this dependence suppressed in the notation.

\subsection{Comparison models and inferential target}
\label{subsec:covshift-arrival-model}

Let \(\vect{\theta}=(\theta_1,\ldots,\theta_K)^\top\) collect the human-reference scores of the historical models. A historical human comparison follows the BTL model
\begin{equation}\label{eq:human-model-new}
\Prob(Y_r^{\Hum}=1\mid I_r^{\Hum}=i,J_r^{\Hum}=j)
=
\sigmoid(\theta_i-\theta_j),
\end{equation}
where \(\sigmoid(z)=(1+e^{-z})^{-1}\). We use \(\vect1_d\), \(\vect0_d\), and \(\mat I_d\) for the all-ones vector, zero vector, and identity matrix of dimension \(d\), respectively. Because the comparison probabilities depend only on score differences, we impose the normalization \(\vect{1}_K^\top\vect{\theta}=0\).

The inferential target \(\theta_\star\) is the new model's score under the human BTL model. If a human comparison between \(\star\) and a historical model \(j\) were observed, its outcome would satisfy
\(\Prob(Y^{\Hum,\star}=1\mid J=j)=\sigmoid(\theta_\star-\theta_j)\).
The new model's position among the historical models is determined by the contrasts \(\theta_\star-\theta_j\), \(j\in\mathcal I\).

For judge \(m\), a historical comparison follows
\begin{equation}\label{eq:hist-surrogate-model-new}
\Prob(Y_r^{C_m}=1\mid \vect{X}_{I_r^m,r}^{m}=\vect{x}_i,\vect{X}_{J_r^m,r}^{m}=\vect{x}_j,
I_r^m=i,J_r^m=j)
=
\sigmoid\!\left[
 c_m(\theta_i-\theta_j)+h_m(\vect{x}_i,\vect{x}_j)
\right].
\end{equation}
Here \(c_m>0\) measures judge \(m\)'s sensitivity to the human score difference, while \(h_m\) captures its feature-dependent bias. We represent \(h_m\) as
\(h_m(\vect{x},\vect{x}')=\psi_m^*(\vect{x})-\psi_m^*(\vect{x}')\),
for \(\vect x,\vect x'\in\mathcal X\), where \(\psi_m^*:\mathcal X\to\mathbb R\) is a bounded continuous single-response bias function for judge \(m\).
The comparison probabilities depend on \(\psi_m^*\) only through the bias contrast \(h_m\), so adding a constant to \(\psi_m^*\) leaves them unchanged. We restrict \(h_m\) to this difference form so that reversing the two responses reverses its sign and the same function can be evaluated on new-model responses. This form assumes that comparison bias is the difference between two single-response contributions and excludes pairwise interactions that cannot be represented in this way. Residual item-specific effects not represented by the observed features are also excluded.

We assume that the same \((c_m,h_m)\) applies to the historical and new-model comparisons evaluated by judge \(m\), although their feature distributions may differ.
The corresponding model for a new-model comparison is
\begin{equation}\label{eq:arrival-surrogate-model-new}
\Prob(Y_s^A=1\mid \vect{X}_{\star,s}^{\star}=\vect{x}_\star,
\vect{X}_{J_s,s}^{\star}=\vect{x}_j,J_s=j,G_s=m)
=
\sigmoid\!\left[
 c_m(\theta_\star-\theta_j)+h_m(\vect{x}_\star,\vect{x}_j)
\right].
\end{equation}
For vectors, \(\|\cdot\|_2\) and \(\|\cdot\|_\infty\) denote the Euclidean and maximum norms; for functions, \(\|\cdot\|_\infty\) is the supremum norm. For fixed \(B_\theta>0\), \(\underline\theta_\star<\overline\theta_\star\), and \(0<\underline c<\overline c<\infty\), we require \(\vect\theta\in\Theta_0\), \(\theta_\star\in\Theta_\star\), and \(\vect c=(c_1,\ldots,c_M)^\top\in\mathcal C\), where
\[
\begin{gathered}
\Theta_0:=\{\vect{t}\in\mathbb R^K:\vect{1}_K^\top \vect{t}=0,\ \|\vect{t}\|_\infty\leq B_\theta\},\ \
\Theta_\star=[\underline\theta_\star,\overline\theta_\star],
\quad
\mathcal C=[\underline c,\overline c]^M.
\end{gathered}
\]

For the feature-dependent biases, write \(\vect h=(h_1,\ldots,h_M)\) and let \(C(\mathcal X)\) denote the continuous real-valued functions on \(\mathcal X\). Given a fixed \(0<B_\psi<\infty\), we define the function classes \(\mathcal H_m\) by
\begin{equation*}
 \mathcal H_m
 :=\left\{u:
 u(\vect{x},\vect{x}')=\psi(\vect{x})-\psi(\vect{x}')
 \text{ for }\psi\in C(\mathcal X),\ \|\psi\|_\infty\leq B_\psi\right\},
 \qquad
 h_m\in\mathcal H_m.
\end{equation*}
The full parameter space is
\(\mathcal P:=\Theta_0\times\mathcal C\times\prod_{m=1}^M\mathcal H_m\times\Theta_\star\),
with \((\vect{\theta},\vect{c},\vect{h},\theta_\star)\in\mathcal P\).

\subsection{Identification}
\label{sec:identification-v2}

Human comparisons anchor the historical scores, but judge sensitivities can be identified only if human score differences cannot be absorbed into feature-dependent biases. Identifying the new model's score further requires the bias functions learned from historical comparisons to be determined on its feature distribution.

For a generic historical comparison from judge \(m\), write \((I,J)\) for its model indices and \(\vect{Z}_{C_m}=(\vect{X}_{I}^{m},\vect{X}_{J}^{m})\) for its pair features. For a generic new-model comparison conditional on \(G=m\), write \(\vect{Z}_{A_m}=(\vect{X}_{\star}^{\star},\vect{X}_{J}^{\star})\). Let \(P_{C_m}^{Z}\) and \(P_{A_m}^{Z}\) denote the distributions of the pair features under \(P_{C_m}\) and \(P_{A_m}\), respectively. For a square-integrable function \(f\) of the full data from a comparison, define
\(\|f\|_{C_m}^2:=\E_{C_m}(f^2)\) and
\(\|f\|_{A_m}^2:=\E_{A_m}(f^2)\).
For a function of the pair features alone, these are the \(L_2\) norms under \(P_{C_m}^{Z}\) and \(P_{A_m}^{Z}\).
Let \(\mathcal G_{\Hum}\) be the undirected graph on \(\mathcal I\) whose edges are the unordered human comparison pairs with positive probability under \(P_H\).

\begin{assumption}
\label{ass:historical-anchor-v2}
The graph \(\mathcal G_{\Hum}\) is connected.
\end{assumption}

Under the zero-sum normalization, Assumption~\ref{ass:historical-anchor-v2} identifies the historical reference scores from human comparisons \citep{negahban2012iterative,chen2015spectral,chen2019spectral}. We retain this assumption for the estimator and inference theory below. It is not necessary for identification from the combined human and judge comparisons: Supplementary Theorem~C.1 gives sufficient conditions that identify the historical parameters and the new-model score even when human comparisons involve only one reference pair with distinct scores.

\begin{assumption}
\label{ass:calibration-separation-v2}
For each judge \(m\), there exists a constant \(\kappa_{c,m}>0\) such that,
for every \(a_m\in\mathbb R\) and \(\varphi_m\in C(\mathcal X)\),
\begin{equation*}
 \left\|a_m(\theta_I-\theta_J)
       +\varphi_m(\vect{X}_{I}^{m})-\varphi_m(\vect{X}_{J}^{m})\right\|_{C_m}
 \geq \kappa_{c,m}|a_m|.
\end{equation*}
\end{assumption}

Assumption~\ref{ass:calibration-separation-v2} separates the historical human score difference from feature-dependent bias contrasts by at least \(\kappa_{c,m}>0\) in \(L_2(P_{C_m})\). Once \(\vect\theta\) is identified, it distinguishes a change in judge sensitivity \(c_m\) from a change in the feature-dependent bias contrast \(h_m\), thereby identifying both components on the historical support.

\begin{assumption}
\label{ass:arrival-overlap-v2}
There exists a constant \(\overline w<\infty\) such that, for every \(m\in\mathcal M_A\),
\begin{equation*}
 P_{A_m}^{Z}\ll P_{C_m}^{Z},\qquad
 \frac{dP_{A_m}^{Z}}{dP_{C_m}^{Z}}(\vect{z})
 \leq\overline w
 \quad P_{C_m}^{Z}\text{-a.e. }\vect{z}.
\end{equation*}
Here \(\ll\) denotes absolute continuity.
\end{assumption}

Absolute continuity ensures that the bias contrast \(h_m\) identified under the historical feature law is also determined under the new-model feature law. The bounded density ratio additionally transfers \(L_2\) error bounds from the historical to the new-model feature law, as used in the estimation and inference results.

\begin{theorem}[Historical and arrival identification]
\label{thm:historical-arrival-identification-v2}
Under Assumptions~\ref{ass:historical-anchor-v2},
\ref{ass:calibration-separation-v2}, and \ref{ass:arrival-overlap-v2},
consider two parameter collections
\(\bigl(\vect{\theta},\{c_m,h_m\}_{m=1}^M,\theta_\star\bigr)\) and
\(\bigl(\bar{\vect{\theta}},\{\bar c_m,\bar h_m\}_{m=1}^M,\bar\theta_\star\bigr)\)
that satisfy the comparison models \eqref{eq:human-model-new}, \eqref{eq:hist-surrogate-model-new}, and \eqref{eq:arrival-surrogate-model-new} and the restrictions
\(\vect{\theta},\bar{\vect{\theta}}\in\Theta_0\),
\(c_m,\bar c_m>0\), and \(h_m,\bar h_m\in\mathcal H_m\) for every \(m\).
If they induce the same joint distribution of
\(\bigl(\mathcal O_{\Hum},\{\mathcal O_{C_m}\}_{m=1}^M,
\mathcal O_A\bigr)\), then
\(\bar{\vect{\theta}}=\vect{\theta}\), \(\bar c_m=c_m\), and
\(\bar h_m(\vect z)=h_m(\vect z)\)
\(P_{C_m}^{Z}\)-almost surely for every \(m\), while
\(\bar\theta_\star=\theta_\star\) and
\(\bar h_m(\vect z)=h_m(\vect z)\)
\(P_{A_m}^{Z}\)-almost surely for every \(m\in\mathcal M_A\).
\end{theorem}

Theorem~\ref{thm:historical-arrival-identification-v2} establishes that all score contrasts \(\theta_\star-\theta_j\), \(j\in\mathcal I\), and hence the new model's position among the historical models, are identified. Under the stated assumptions, this requires no human comparisons involving the new model and allows different feature distributions in the historical and new-model LLM-judge comparisons.

\section{Estimation and inference}
\label{sec:estimation-v2}
\label{sec:method}

\subsection{Two-step estimator}
\label{subsec:estimation-v2}

Since a joint likelihood fit to both historical and new-model
comparisons would re-estimate the historical scores and judge parameters
for every new model, we use a two-step estimator that reuses the historical
fit across evaluations.
The first step estimates the historical scores and judge
parameters from the historical human and LLM comparisons. The second step holds
these estimates fixed and estimates \(\theta_\star\) from the new model's
LLM-judge comparisons through a one-dimensional optimization, reducing the computation required for each new model.

We estimate \(\vect\theta\) from the human comparisons \(\mathcal O_{\Hum}\)
under the BTL model in \eqref{eq:human-model-new} by maximum likelihood,
using the logistic loss
\mbox{\(\mathcal L_{\rm log}(y,\eta)=-y\eta+\log(1+e^\eta)\)}:
\begin{equation}
\widetilde{\vect{\theta}}
 \in
 \argmin_{\vect{\theta}\in\Theta_0}
\sum_{r=1}^{n_{\Hum}}
\mathcal L_{\rm log}\!\left\{Y_r^{\Hum},
 \theta_{I_r^{\Hum}}-\theta_{J_r^{\Hum}}\right\}.
\label{eq:human-estimation-v2}
\end{equation}

With \(\widetilde{\vect\theta}\) fixed, we estimate \(c_m\) and \(h_m\)
from judge \(m\)'s historical comparisons \(\mathcal O_{C_m}\).
Using the difference representation of \(h_m\) in Section~\ref{subsec:covshift-arrival-model}, we estimate its single-response bias function \(\psi_m^*\) with a deep neural network (DNN) using rectified linear unit (ReLU) activations.
The neural network function \(\psi\) is constrained to the bounded sparse network class
\(\mathcal F_m^{\rm NN}\). The network architecture and sample-size-dependent
choices of depth, width, and sparsity are detailed in
Supplementary Section~C.1.
Given \(\widetilde{\vect{\theta}}\), we define the judge-specific loss as
\begin{equation*}
 \mathcal L_{C_m}(c,\psi;\widetilde{\vect{\theta}})
 =\sum_{r=1}^{n_{C_m}}
 \mathcal L_{\rm log}\!\left[
 Y_r^{C_m},\,
 c\{\widetilde\theta_{I_r^m}-\widetilde\theta_{J_r^m}\}
 +\psi(\vect{X}_{I_r^m,r}^m)-\psi(\vect{X}_{J_r^m,r}^m)
 \right]
\end{equation*}
and minimize it:
\begin{equation}
 (\widetilde c_m,\widetilde\psi_m)
 \in\argmin_{\substack{c\in[\underline c,\overline c]\\
                       \psi\in\mathcal F_m^{\rm NN}}}
 \mathcal L_{C_m}(c,\psi;\widetilde{\vect{\theta}}).
\label{eq:calibration-estimation-v2}
\end{equation}
Set \(\widetilde h_m(\vect{x},\vect{x}')
:=\widetilde\psi_m(\vect{x})-\widetilde\psi_m(\vect{x}')\).
Let
\(\widetilde{\vect c}=(\widetilde c_1,\ldots,\widetilde c_M)^\top\) and
\(\widetilde{\vect h}=(\widetilde h_1,\ldots,\widetilde h_M)^\top\).

We use this historical fit because fixing \(\widetilde{\vect\theta}\)
allows the judge models to be fitted independently, simplifying the optimization.
The orthogonal correction in Section~\ref{sec:targeted-inference-v2} then combines
information from all three samples and yields a semiparametrically efficient
estimator under the conditions of Theorem~\ref{thm:human-scale-information-limit-v2}.
Joint minimization of the historical human and judge losses is also possible;
Supplementary Section~C.4 gives the objective,
and Section~A.4 compares the two approaches.

We plug \((\widetilde{\vect\theta},\widetilde{\vect c},\widetilde{\vect h})\) into the new-model comparison likelihood and estimate \(\theta_\star\) from all \(n_\star\) comparisons:
\begin{equation}
\widetilde\theta_\star
\in\argmin_{t\in\Theta_\star}
\sum_{s=1}^{n_\star}
\mathcal L_{\rm log}\!\left[
Y_s^A,\,
\widetilde c_{G_s}\{t-\widetilde\theta_{J_s}\}
+\widetilde h_{G_s}(\vect{X}_{\star,s}^{\star},\vect{X}_{J_s,s}^{\star})
\right].
\label{eq:arrival-estimation-v2}
\end{equation}

\subsection{Orthogonal correction and inference}
\label{sec:targeted-inference-v2}

Our goal is inference on the new model's human-reference score
\(\theta_\star\) or its difference \(\theta_\star-\theta_j\) from the
score of a prespecified historical model \(j\). Write
\(\vect\theta=\mat B\vect\vartheta\) to enforce the constraint
\(\vect1_K^\top\vect\theta=0\), where the columns of \(\mat B\)
form an orthonormal basis of the zero-sum subspace. Set \(p:=K+M\),
\(\vect\beta=(\theta_\star,\vect\vartheta^\top,\vect c^\top)^\top\in\mathbb R^p\),
and \(T_\ell=\vect\ell^\top\vect\beta\).
The choices
\((1,\vect0_{K-1}^\top,\vect0_M^\top)^\top\) and
\((1,-\mat B_{j,\cdot},\vect0_M^\top)^\top\) select the score and
contrast, respectively. Here \(\mat B_{j,\cdot}\) is row \(j\) of \(\mat B\).
{For the inference of \(\theta_\star\), the historical scores, judge sensitivities,
and bias functions are nuisance quantities. Their estimation errors propagate to the
initial estimate in \eqref{eq:arrival-estimation-v2}.}
{In this subsection, we correct these excessive errors using weighted residuals from all three samples,
with weights chosen jointly to remove the first-order effects of these nuisance
estimation errors.}

\begingroup

To implement this construction, we use cross-fitting to evaluate each correction on observations excluded
from fitting the models and weights \citep{chernozhukov2018double}.
Split each source \(q\in\mathcal S\) into \(K_{\cf}\geq2\) equal folds
\(\mathcal F_{q,k}\), with fixed \(K_{\cf}\) dividing every \(n_q\).
Fit the two-step estimator without fold \(k\), denoting the fits by
\((\widetilde{\vect\beta}^{(-k)},\widetilde{\vect h}^{(-k)})\).
Let \(\eta_q(\vect\beta,\vect h)\) be the logit in source \(q\)'s model
\eqref{eq:human-model-new}--\eqref{eq:arrival-surrogate-model-new}.
The fitted logit \(\widehat\eta_{q,r}^{(-k)}\) gives the residual
\(\widehat\varepsilon_{q,r}^{(-k)}:=Y_r^q-\sigmoid(\widehat\eta_{q,r}^{(-k)})\).
Our estimator adds a weighted residual correction to the average initial estimate:
\begin{equation}
 \widehat T_\ell
 =\frac1{K_{\cf}}\sum_{k=1}^{K_{\cf}}\vect\ell^\top\widetilde{\vect\beta}^{(-k)}
 +\sum_{k=1}^{K_{\cf}}\sum_{q\in\mathcal S}
  \sum_{r\in\mathcal F_{q,k}}
 \widehat a_{\ell,q,r}^{(-k)}
  \widehat\varepsilon_{q,r}^{(-k)}.
  \label{eq:feasible-anchor-estimator-v2}
\end{equation}
The key debiasing term
\(\widehat a_{\ell,q,r}^{(-k)}\widehat\varepsilon_{q,r}^{(-k)}\)
combines an outcome residual $\widehat\varepsilon_{q,r}^{(-k)}$ with a weight $\widehat a_{\ell,q,r}^{(-k)}$ that converts it into a correction
for \(T_\ell\). The residual is the observed comparison outcome minus its
fitted win probability. We choose the weights jointly across all three samples
so that the weighted residual correction offsets the first-order error in the
initial target estimate. The specific construction procedures of $\widehat a_{\ell,q,r}^{(-k)}$ are introduced below.

To compute the weights, let \(\widehat{\vect\xi}_{q,r}^{(-k)}\) be the
gradient \(\nabla_{\vect\beta}\eta_{q,r}\) evaluated at the fitted parameters,
and set \(\widehat V_{q,r}^{(-k)}:=\sigmoid'(\widehat\eta_{q,r}^{(-k)})\).
Let \(\vect Z_{q,r}\) denote the ordered pair of left and right response
feature vectors in judge comparison \(r\) from source \(q\).
For functions \(\vect\gamma=(\gamma_1,\ldots,\gamma_M)\) with
\(\gamma_m(\vect x,\vect x')=\varphi_m(\vect x)-\varphi_m(\vect x')\),
where \(\varphi_m\) is a single-response function, define
\begin{equation*}
 D_{q,r}[\vect\gamma]:=\gamma_m(\vect Z_{q,r}),\quad q\in\{C_m,A_m\},
 \qquad D_{H,r}[\vect\gamma]:=0.
\end{equation*}
For judge comparisons, this is the logit change induced by a bias perturbation
in direction \(\gamma_m\). The human model has no judge bias term.
We define the correction weights by
\begin{equation}
 \widehat a_{\ell,q,r}^{(-k)}
 :=\widehat{\vect w}_{\ell,k}^{\top}\widehat{\vect\xi}_{q,r}^{(-k)}
   -D_{q,r}[\widehat{\vect\gamma}_{\ell,k}],
 \label{eq:estimated-riesz-weights-v2}
\end{equation}
where the coefficients and functions are fitted jointly on the training observations:
\begin{equation}
 (\widehat{\vect w}_{\ell,k},\widehat{\vect\gamma}_{\ell,k})
 \in\argmin_{\substack{\vect w\in\mathbb R^p\\
        \gamma_m\in n_{\Hum}^{-1}\Gamma_m^R}}
 \bigg\{
 \frac1{K_{\cf}-1}\sum_{q\in\mathcal S}
  \sum_{r\notin\mathcal F_{q,k}}
  \{\vect w^\top\widehat{\vect\xi}_{q,r}^{(-k)}
    -D_{q,r}[\vect\gamma]\}^2\widehat V_{q,r}^{(-k)}
  -\frac{2}{K_{\cf}}\vect\ell^\top\vect w
 \bigg\}.
 \label{eq:fold-direct-riesz-v2}
\end{equation}
Here \(\Gamma_m^R\) consists of pairwise differences of bounded neural-network
functions, specified in Supplementary Section~C.1.

To gain insight into the form of the debiasing correction, consider the
population weights. The key identities are given below, with
the full optimization problem, function-space definitions, and derivation
deferred to the supplement to keep the main exposition concise.
Write \(\eta_q=\eta_q(\vect\beta,\vect h)\),
\(\vect\xi_q=\nabla_{\vect\beta}\eta_q\), and \(V_q=\sigmoid'(\eta_q)\)
at the truth, and let \(D_q\) denote the same evaluation as \(D_{q,r}\)
at a generic observation from source \(q\).
The true residual is \(\varepsilon_q:=Y^q-\sigmoid(\eta_q)\),
with \(\varepsilon_{q,r}\) denoting its value at comparison \(r\).
The population weights
\(a_{\ell,q}=\vect w_\ell^\top\vect\xi_q-D_q[\vect\gamma_\ell]\)
solve the population criterion \textup{(C.10)} and satisfy
\begin{equation}
 \sum_q n_q\E_q(a_{\ell,q}V_q\vect\xi_q)=\vect\ell,
 \qquad
 \sum_q n_q\E_q\{a_{\ell,q}V_qD_q[\vect g]\}=0
 \quad\text{for all admissible }\vect g.
 \label{eq:split-cancellation-conditions-v2}
\end{equation}
The admissible bias directions \(g_m\) belong to the spaces \(\mathcal T_m\)
generated by continuous pairwise differences, defined in Supplementary
Section~C.3.
The first identity removes the first-order target error from estimating
\(\vect\beta\). The second removes the first-order effect of estimating \(\vect h\).
For each judge, the historical and new-model samples contribute jointly to
this bias correction because they share the same bias function.
Lemma~C.1 gives the derivation and
shows that these weights minimize the first-order variance
\(v_\ell^2:=\sum_q n_q\E_q(a_{\ell,q}^2V_q)\)
among weights satisfying \eqref{eq:split-cancellation-conditions-v2}.
Equivalently, we can also view these as Riesz weights for the target derivative
\citep{chernozhukov2022automatic,chernozhukov2022riesz}.

Based on the debiased estimator in \eqref{eq:feasible-anchor-estimator-v2},
we estimate \(v_\ell^2\) by summing the squared debiasing contributions
from the held-out observations:
\begin{equation}
  \widehat v_\ell^2
 =\sum_{k=1}^{K_{\cf}}\sum_{q\in\mathcal S}
   \sum_{r\in\mathcal F_{q,k}}
   \{\widehat a_{\ell,q,r}^{(-k)}
    \widehat\varepsilon_{q,r}^{(-k)}\}^2.
  \label{eq:feasible-anchor-variance-v2}
\end{equation}
Algorithm~1 in Supplementary
Section~C.2 summarizes the complete
cross-fitted estimation and inference procedure.
\endgroup

\section{Theoretical analysis}
\label{sec:theory-v2}
This section establishes the theoretical guarantees for the proposed estimation and inference procedures. Section~\ref{subsec:estimation-theory-v2} derives convergence rates for the two-step estimator. Section~\ref{subsec:inference-v2} establishes asymptotic normality and consistent variance estimation for the cross-fitted estimator. Section~\ref{subsec:human-scale-lower-bound-v2} establishes its semiparametric efficiency.
\subsection{Rates for the two-step estimator}
\label{subsec:estimation-theory-v2}

\begin{assumption}
\label{ass:estimation-geometry-v2}
For each judge \(m\), there exists \(s_m>0\) such that
\(\psi_m^*\) belongs to the H\"older ball
\(C_{d_x}^{s_m}([-1,1]^{d_x};B_\psi)\) defined in
Supplementary Definition~C.1.
\end{assumption}

Define
\mbox{\(r_{{\rm NN},m}:=\{(\log n_{C_m})/n_{C_m}\}^{1/2}
+n_{C_m}^{-s_m/(2s_m+d_x)}(\log n_{C_m})^2\)}
for \(m=1,\ldots,M\), and
\mbox{\(r_{\Hist}:=\max_{1\leq m\leq M}(1+\kappa_{c,m}^{-1})
\{n_{\Hum}^{-1/2}+r_{{\rm NN},m}\}\)}.
The exponent \(s_m/(2s_m+d_x)\) reflects the difference
structure \(h_m(\vect x,\vect x')=\psi_m^*(\vect x)-\psi_m^*(\vect x')\):
the network approximates a \(d_x\)-dimensional function, although each
comparison contains two feature vectors. Up to logarithmic factors,
\(r_{{\rm NN},m}\) has the standard H\"older rate
\citep{schmidthieber2020relu,farrell2021deep}.
In \(r_{\Hist}\), the term \(n_{\Hum}^{-1/2}\) arises because the estimated
historical scores enter the fit of \(c_m\) and \(h_m\).
The factor \(1+\kappa_{c,m}^{-1}\) increases the bound when changes in
judge sensitivity are difficult to distinguish from changes in
feature-dependent bias.

\begin{theorem}[Rates for the two-step estimator]
\label{thm:estimation-rates-v2}
Suppose \(K\) and \(M\)
are fixed and
\(n_{\Hum}\wedge n_\star\wedge\min_{1\leq m\leq M}n_{C_m}\to\infty\).
Suppose also that Assumptions~\ref{ass:historical-anchor-v2},
\ref{ass:calibration-separation-v2}, \ref{ass:arrival-overlap-v2}, and
\ref{ass:estimation-geometry-v2} hold and the two-step estimator uses the
network class in \textup{(C.4)}. Write
\(\Delta\widetilde c_m:=\widetilde c_m-c_m\) and
\(\Delta\widetilde h_{C_m}:=(\widetilde h_m-h_m)(\vect Z_{C_m})\). For
\(m\in\mathcal M_A\), also write
\(\Delta\widetilde h_{A_m}:=(\widetilde h_m-h_m)(\vect Z_{A_m})\). Then
\begin{align*}
 \|\widetilde{\vect{\theta}}-\vect{\theta}\|_2
 &=O_p(n_{\Hum}^{-1/2}),
 &
 |\widetilde\theta_\star-\theta_\star|
 &=O_p(n_\star^{-1/2}+r_{\Hist}),\\
 \max_{1\leq m\leq M}\!\left\{|\Delta\widetilde c_m|
 +\|\Delta\widetilde h_{C_m}\|_{C_m}\right\}
 &=O_p(r_{\Hist}),
 &
 \max_{m\in\mathcal M_A}\|\Delta\widetilde h_{A_m}\|_{A_m}
 &=O_p(r_{\Hist}).
\end{align*}
\end{theorem}
Although \(\theta_\star\) is scalar, its two-step estimator
inherits the error from fitting the historical scores, judge sensitivities,
and bias functions. The bound gives an \(n_\star^{-1/2}\) rate when
\(r_{\Hist}=O(n_\star^{-1/2})\). When historical fitting dominates the bound,
increasing the new-model sample size alone reduces only the smaller
sampling term. The bounded density ratio in
Assumption~\ref{ass:arrival-overlap-v2} ensures that the change in feature
distribution does not further reduce the convergence rate.

\subsection{Asymptotic normality and variance estimation}
\label{subsec:inference-v2}

We establish asymptotic normality and consistent variance estimation for the cross-fitted estimator \(\widehat T_\ell\) constructed in Section~\ref{sec:targeted-inference-v2}.
Let \(\mathcal S_J:=\mathcal S\setminus\{H\}\) be the set of LLM judge
sources, with total sample size
\(N_J:=\sum_{m=1}^M n_{C_m}+\sum_{m\in\mathcal M_A}n_{A_m}\) and source
shares \(\pi_q:=n_q/N_J\) for \(q\in\mathcal S_J\).
Throughout the remainder of this section, assume
\begin{equation}
 n_{\Hum}\to\infty,
 \qquad \liminf_n N_J/n_{\Hum}>0,
 \qquad \min_{q\in\mathcal S_J}\pi_q
 \geq\underline\pi.
 \label{eq:inference-sample-size-conditions-v2}
\end{equation}
Here \(\underline\pi>0\) is a fixed constant, and the inequalities hold for all sufficiently large \(n\). We require the human sample size \(n_{\Hum}\) to diverge and the total judge sample size \(N_J\) to grow at least as fast as \(n_{\Hum}\). The bound \(\pi_q\geq\underline\pi\) ensures that each source contributes a nonvanishing fraction of comparisons.

\begingroup

To characterize the estimator's variance, we first distinguish the information
supplied by human and LLM-judge comparisons.
LLM comparison probabilities remain unchanged if all human-reference scores
are multiplied by a common positive factor and judge sensitivities are divided
by that factor. The corresponding scale direction and its unit normalization are
\[
\vect s_0:=(\theta_\star,\vect\vartheta^\top,-\vect c^\top)^\top,
\qquad
\vect v_0:=\vect s_0/\|\vect s_0\|_2.
\]
Only human comparisons provide information along \(\vect v_0\), and
\(\vect\ell^\top\vect v_0\) measures the target's sensitivity to this rescaling. To quantify the LLM-judge information beyond this scale direction, we define the
following information matrix per comparison, accounting for the unknown bias functions:
\begin{equation*}
 \vect b^\top\mat{\mathcal J}_J\vect b
 :=\inf_{\vect\gamma\in\prod_{m=1}^M\mathcal T_m}
 \sum_{q\in\mathcal S_J}\pi_q\E_q\!\left[
 V_q\{\vect\xi_q^\top\vect b-D_q[\vect\gamma]\}^2\right].
\end{equation*}
The same bias perturbation is used for each judge's historical and new-model
samples. This matrix satisfies \(\mat{\mathcal J}_J\vect v_0=\vect0\).
The following assumption requires its quadratic form to be uniformly positive
over unit directions orthogonal to \(\vect v_0\).
\endgroup

\begin{assumption}
\label{ass:joint-profile-information-v2}
There is a constant \(\kappa_J>0\) such that, for all sufficiently large
\(n\), \(\vect b^\top\mat{\mathcal J}_J\vect b\geq\kappa_J\) whenever
\(\|\vect b\|_2=1\) and \(\vect b^\top\vect v_0=0\).
\end{assumption}

Let \((\vect w_\ell,\vect\gamma_\ell)\) be a population Riesz solution
to \textup{(C.10)}, and set
\(\lambda_{\ell,m}:=n_{\Hum}\gamma_{\ell,m}\).
Human comparisons supply the information in the common scale direction, so
we use the rescaled functions \(\lambda_{\ell,m}\) to retain this contribution
even when \(N_J/n_{\Hum}\) diverges.

\begin{assumption}
\label{ass:riesz-dnn-approximation-v2}
There is a fixed \(B_\lambda<\infty\) and, for each \(m\), a fixed
\(t_m>d_x/2\) such that, for every sufficiently large \(n\), there is a
continuous function \(\lambda_{\ell,m}^\circ\) satisfying
\begin{equation*}
 \lambda_{\ell,m}(\vect x,\vect x')
 =\lambda_{\ell,m}^\circ(\vect x)-\lambda_{\ell,m}^\circ(\vect x'),
 \qquad
 \lambda_{\ell,m}^\circ\in C_{d_x}^{t_m}([-1,1]^{d_x};B_\lambda).
\end{equation*}
The equality holds \(P_{C_m}^Z\)-almost surely and, for
\(m\in\mathcal M_A\), \(P_{A_m}^Z\)-almost surely.
\end{assumption}

For judge \(m\), estimating \(\lambda_{\ell,m}\) uses its historical
and available new-model comparisons. Write their combined sample size as
\(n_{{\rm Rz},m}:=n_{C_m}+1\{m\in\mathcal M_A\}n_{A_m}\), and set
\[
 r_\eta:=n_{\Hum}^{-1/2}+N_J^{-1/2}+\max_m r_{{\rm NN},m},
 \qquad
 r_\lambda:=\max_m n_{{\rm Rz},m}^{-t_m/(2t_m+d_x)}
 (\log n_{{\rm Rz},m})^2.
\]
\begingroup

The rate \(r_\eta\) has the same order as
\(n_\star^{-1/2}+r_{\Hist}\), the error bound for the initial new-model score
in Theorem~\ref{thm:estimation-rates-v2}, since
\(n_\star\asymp N_J\) under the source-share condition
\eqref{eq:inference-sample-size-conditions-v2}.
It summarizes the error rate for the initial estimates of the scores, judge
sensitivities, and bias functions. In the debiasing step, \(r_\lambda\) is the
network error component for estimating the functions \(\lambda_{\ell,m}\)
used to construct the Riesz weights.
\endgroup
To make the remainder negligible relative to \(v_\ell\), we require
\begin{equation}
 \Delta_{\ell,n}:=
 \frac{N_J}{n_{\Hum}v_\ell}\,
 r_\eta(r_\eta+r_\lambda)\longrightarrow0.
 \label{eq:main-product-rate-v2}
\end{equation}
{ The product \(r_\eta(r_\eta+r_\lambda)\) contains the squared initial
estimation error and its interaction with the network error term \(r_\lambda\).} The factor
\(N_J/n_{\Hum}\) arises from summing the estimation-error contributions
over \(N_J\) judge comparisons, with Riesz weights bounded at order
\(n_{\Hum}^{-1}\).

\begin{theorem}[Asymptotic normality and variance estimation]
\label{thm:attainment-v2}
Let \(\vect\ell\in\mathbb R^{K+M}\setminus\{\vect0\}\) be fixed.
Suppose \eqref{eq:inference-sample-size-conditions-v2}, the conditions of
Theorem~\ref{thm:estimation-rates-v2}, and Assumptions \ref{ass:joint-profile-information-v2}--\ref{ass:riesz-dnn-approximation-v2}
hold, the networks used to estimate \(\lambda_{\ell,m}^{\circ}\) have the
architecture in \textup{(C.6)}, and
\(\Delta_{\ell,n}\to0\).
Then
\begin{equation*}
 v_\ell^{-1}\left\{
 \widehat T_\ell-\vect\ell^\top\vect\beta
 -\sum_{q\in\mathcal S}\sum_{r=1}^{n_q}
 a_{\ell,q,r}\varepsilon_{q,r}\right\}=o_p(1).
\end{equation*}
Moreover,
\begin{equation*}
 \frac{\widehat T_\ell-\vect\ell^\top\vect\beta}{v_\ell}
 \rightsquigarrow N(0,1),
 \qquad \frac{\widehat v_\ell^2}{v_\ell^2}\longrightarrow_p1,
 \qquad
 \frac{\widehat T_\ell-\vect\ell^\top\vect\beta}{\widehat v_\ell}
 \rightsquigarrow N(0,1).
\end{equation*}
\end{theorem}

Let \(\mat{\mathcal J}_H\) denote the information
per human comparison and \(\mat{\mathcal J}\) the total information after
accounting for the unknown judge bias functions:
$
 \mat{\mathcal J}_H:=\E_H(V_H\vect\xi_H\vect\xi_H^\top),
 \mat{\mathcal J}:=n_{\Hum}\mat{\mathcal J}_H+N_J\mat{\mathcal J}_J.
$
Under the conditions of the theorem, \(\mat{\mathcal J}\) is positive
definite, and the variance satisfies
\begin{equation}
 v_\ell^2=\vect\ell^\top\mat{\mathcal J}^{-1}\vect\ell
 \asymp
 \frac{(\vect\ell^\top\vect v_0)^2}{n_{\Hum}}
 +\frac{\|(\mat I_p-\vect v_0\vect v_0^\top)\vect\ell\|_2^2}{N_J}.
 \label{eq:target-variance-orders-v2}
\end{equation}
For a fixed target with \(\vect\ell^\top\vect v_0\ne0\), the variance is of order \(n_{\Hum}^{-1}\).
When \(\vect\ell^\top\vect v_0=0\), the target is locally unaffected by
rescaling, and its variance is of order \(N_J^{-1}\).

The case \(\vect\ell^\top\vect v_0=0\) arises naturally when comparing
models with equal scores. For a
fixed historical model \(j\), consider \(T_\ell=\theta_\star-\theta_j\),
for which
\(\vect\ell^\top\vect v_0=(\theta_\star-\theta_j)/\|\vect s_0\|_2\).
When \(\theta_\star=\theta_j\), rescaling preserves the equality, so
this derivative vanishes and \(v_\ell^2\asymp N_J^{-1}\).
{ The asymptotic normality in Theorem~\ref{thm:attainment-v2} applies both when
\(N_J\asymp n_{\Hum}\) and when \(n_{\Hum}=o(N_J)\). The orthogonal correction removes the first-order effects
of initial estimation errors, while condition~\eqref{eq:main-product-rate-v2}
requires the resulting remainder to be negligible relative to the target's
standard deviation \(v_\ell\). Thus, the variance orders in
\eqref{eq:target-variance-orders-v2} determine how small the errors in the
initial estimates and fitted Riesz weights must be in each case.}

When \(N_J\asymp n_{\Hum}\), the remainder condition
is satisfied for both types of target if \(r_\eta\) and \(r_\lambda\)
are \(o(n_{\Hum}^{-1/4})\). Since \(\pi_q\geq\underline\pi>0\) ensures
that each judge source has a sample size of order \(n_{\Hum}\), these
rates hold when \(s_m>d_x/2\) and
\(t_m>d_x/2\) for every judge.

When \(n_{\Hum}=o(N_J)\), the two standard deviations have different
orders, so the rate requirements depend on the target. To express these
requirements in terms of sample sizes and smoothness, let
\(s_{\min}:=\min_{1\leq m\leq M}\{s_m,t_m\}>d_x/2\).
If \(N_J=O(n_{\Hum}^{1+\kappa})\) for a fixed \(\kappa>0\), then
\(\Delta_{\ell,n}\to0\) whenever
\[
 \kappa<
 \begin{cases}
 \displaystyle
 \min\left\{\frac12,\frac{s_{\min}}{d_x}-\frac12\right\},
 &\vect\ell^\top\vect v_0\ne0,\\[6pt]
 \displaystyle
 \min\left\{\frac13,
       \frac{2s_{\min}-d_x}{2s_{\min}+3d_x}\right\},
 &\vect\ell^\top\vect v_0=0.
 \end{cases}
\]
The second bound is more restrictive because targets with
\(\vect\ell^\top\vect v_0=0\) have the smaller standard deviation
\(N_J^{-1/2}\), and therefore require a smaller remainder. For example,
when \(s_{\min}=d_x\), the bounds allow \(\kappa<1/2\) and
\(\kappa<1/5\), respectively. These are sufficient bounds for the stated
rates, with strict inequalities to accommodate the logarithmic
factors.

We use the contrast inference in Theorem~\ref{thm:attainment-v2} to
construct a confidence set for the new model's rank. For \(j=1,\ldots,K\),
let \(\Delta_j=\theta_\star-\theta_j\) and set
\(\vect\ell_j=(1,-\mat B_{j,\cdot},\vect0_M^\top)^\top\), so that
\(\vect\ell_j^\top\vect\beta=\Delta_j\). Under the theorem's conditions,
\(\widehat\Delta_j:=\widehat T_{\ell_j}\) is asymptotically normal with
consistently estimated variance \(\widehat s_j^2:=\widehat v_{\ell_j}^2\),
which accounts jointly for uncertainty in the new and historical scores.

Using these estimates, form
\(p_j^-=\Phi(\widehat\Delta_j/\widehat s_j)\) for
\(H_j^-:\Delta_j\geq0\) against \(\Delta_j<0\), and
\(p_j^+=1-\Phi(\widehat\Delta_j/\widehat s_j)\) for
\(H_j^+:\Delta_j\leq0\) against \(\Delta_j>0\), where \(\Phi\) is the
standard normal distribution function. Set both $p$-values to one if
\(\widehat s_j=0\). We apply the Holm step-down procedure
\citep{holm1979simple} separately to the \(K\) lower-tail $p$-values
\(\{p_j^-\}_{j=1}^K\) and the \(K\) upper-tail $p$-values
\(\{p_j^+\}_{j=1}^K\), using level \(\alpha/2\) for each family.
Within each family, sort the $p$-values in increasing order and reject
successively while \(p_{(k)}\leq\alpha/\{2(K-k+1)\}\), stopping at the
first failure.
Let \(\mathcal R_-\) and \(\mathcal R_+\) contain the indices \(j\) for
which \(H_j^-\) and \(H_j^+\) are rejected, respectively. These rejections
identify historical models with higher and lower scores than the new
model, respectively, yielding the following lower and upper bounds on
its rank \citep{mogstad2024ranks}:
\begin{equation}
 \mathcal C_{\rm rank}
 =\{1+|\mathcal R_-|,\ldots,K+1-|\mathcal R_+|\}.
 \label{eq:simulation-rank-set}
\end{equation}

\begin{corollary}[Asymptotic rank coverage]
\label{cor:rank-coverage}
Let \(K\) and \(\alpha\in(0,1)\) be fixed. Suppose the conditions of
Theorem~\ref{thm:attainment-v2} hold for each contrast loading \(\vect\ell_j\),
\(j=1,\ldots,K\). With
\(r_\star:=1+\sum_{j=1}^K1\{\theta_j>\theta_\star\}\), it holds that
\begin{equation*}
 \liminf_{n\to\infty}\Pr\{r_\star\in\mathcal C_{\rm rank}\}
 \geq1-\alpha.
\end{equation*}
\end{corollary}

The coverage guarantee allows tied scores and does not require
independence among the contrast estimates.

\subsection{Semiparametric efficiency}
\label{subsec:human-scale-lower-bound-v2}

The model in Section~\ref{sec:model} is semiparametric in the sense that
\(\vect\beta\) contains finitely many score and judge-sensitivity parameters,
whereas the unknown judge bias functions are infinite-dimensional nuisance
parameters. The covariate and pair-sampling distributions are also left
unspecified, subject to the model restrictions.

To assess efficiency, we compare the asymptotic variances of
estimators of \(T_\ell\) while treating these functions and distributions
as unknown. For an asymptotically linear estimator, the leading error is
a sum of mean-zero contributions from individual comparisons, as in
Theorem~\ref{thm:attainment-v2}, so its asymptotic variance can be evaluated
from this sum. { An estimator is called \emph{regular} if, under
every admissible local perturbation of the model, its distribution after
centering at the perturbed target and dividing by \(v_\ell\) converges to
the same limiting distribution as under the original model.
The admissible perturbations and formal definition are given in Supplementary
Section~E.5.1.} The semiparametric efficiency bound \(V_{\ell,\mathrm{eff}}\)
is a lower bound on the asymptotic variance of all the regular and asymptotically
linear estimators of \(T_\ell\) under the semiparametric model.
We will show that, under the conditions below, ANCHOR is regular and its
asymptotic variance \(v_\ell^2\) given in Section~\ref{subsec:inference-v2}
attains this bound.

Note that only human comparisons provide information along the
rescaling direction \(\vect v_0\). Restricting attention to this direction
gives the variance lower bound
\mbox{\(B_\ell:=(\vect\ell^\top\vect v_0)^2/(n_{\Hum}\vect v_0^\top\mat{\mathcal J}_H\vect v_0)\)}.
The numerator is the target's squared sensitivity to
rescaling, and the denominator is the information supplied by human
comparisons in that direction.
The bound \(V_{\ell,\mathrm{eff}}\) also accounts for uncertainty in the other
parameter directions. The following theorem identifies this bound and
quantifies its difference from \(B_\ell\).
Throughout this subsection, efficiency is understood pointwise at a fixed
interior point of the parameter space.

\begin{theorem}[Efficiency bound]
\label{thm:human-scale-information-limit-v2}
Let \(\vect\ell\in\mathbb R^{K+M}\setminus\{\vect0\}\) be fixed.
At a fixed interior point of the parameter space, suppose the conditions of
Theorem~\ref{thm:attainment-v2} hold.
Then, among regular and asymptotically linear estimators of \(T_\ell\),
the semiparametric efficiency bound satisfies
\begin{equation*}
 V_{\ell,\mathrm{eff}}
 =
 v_\ell^2
 =
 \vect\ell^\top\mat{\mathcal J}^{-1}\vect\ell>0,
 \qquad
 0\leq V_{\ell,\mathrm{eff}}-B_\ell
 \leq C_\ell/N_J,
\end{equation*}
for some constant \(C_\ell<\infty\).
Moreover, ANCHOR is regular and attains
\(V_{\ell,\mathrm{eff}}\).
\end{theorem}

The equality \(V_{\ell,\mathrm{eff}}=v_\ell^2\) shows that no estimator in the
stated class has a
smaller asymptotic variance than ANCHOR. The second conclusion bounds the
gap between \(V_{\ell,\mathrm{eff}}\) and \(B_\ell\).
For a fixed target with \(\vect\ell^\top\vect v_0\ne0\),
\(B_\ell\asymp n_{\Hum}^{-1}\), whereas the gap
\(V_{\ell,\mathrm{eff}}-B_\ell\) is at most of order \(N_J^{-1}\).
When \(N_J\asymp n_{\Hum}\), the gap can be of the same order as
\(B_\ell\). When \(n_{\Hum}=o(N_J)\), it becomes negligible relative to
\(B_\ell\), so that
\(V_{\ell,\mathrm{eff}}/B_\ell\longrightarrow1\).
The leading term of the efficiency bound is therefore
determined by human comparisons.  For targets with \(\vect\ell^\top\vect v_0=0\),
\(B_\ell=0\) because the target is locally unaffected by rescaling.
The bound \(V_{\ell,\mathrm{eff}}\) remains positive, with
\(V_{\ell,\mathrm{eff}}\asymp N_J^{-1}\), so the attainable precision
continues to improve with the number of LLM comparisons under the stated
conditions.

\section{Simulation study}\label{sec:simulation}\label{subsec:sim-setup}

We assess the performance of ANCHOR in score estimation, interval coverage, and ranking, and compare it with the following methods. Implementation details are given in Supplementary Section~A.8. We include two benchmarks with additional information: HumanOnly, which uses $150$ additional human comparisons involving the new model under the default setting, and Oracle, which uses the true nuisance functions. JointBTL pools human and judge comparisons under a common BTL model \citep{bradley1952rank}, ignoring judge-specific sensitivities and feature-dependent bias. DirectScore uses historical estimates of judge sensitivities and bias functions in a plug-in estimating equation for the new model's score, without joint orthogonal correction. IPWOnly adds a density-weighted historical refit to examine whether reweighting for feature shift improves plug-in estimation. HistoricalPPI adapts prediction-powered inference \citep{angelopoulos2023ppi} by adding an average residual correction from historical human labels, without transporting it to the arrival feature distribution. The ranking comparisons also include HistoricalPPR and LogisticJudge, which adapt prediction-powered ranking \citep{chatzi2024ppr} and sigmoid recalibration \citep{guo2017calibration}, respectively, to our setup. Additional comparisons with PooledLLM, ScalarJudge, and PerJudgeBTL are reported in Supplementary Section~A.2.

We simulate data from the models in Section~\ref{sec:model} with \(K=10\) and \(M=3\). In each replication, historical scores are drawn from \(N(0,1)\) and centered at zero, and judge sensitivities are drawn from \(\operatorname{Unif}(0.5,2)\). We sample historical pairs and arrival opponents uniformly, assigning historical comparisons to humans with probability \(\nu_{\Hum}=0.3\) and all other comparisons to uniformly sampled judges.

The main designs use linear judge-bias functions of three response features. Model-specific tanh networks generate these features conditional on a shared prompt vector drawn from \(N(0,\mat I_2)\), with independent \(N(0,0.3^2\mat I_3)\) noise. The fitted bias and Riesz contrasts are linear in response-feature differences. We rescale the bias coefficients using \(\rho\) to control bias magnitude relative to the score term. The new model's response features receive an additional mean shift of magnitude \(\tau\sqrt{5}\).

Unless varied below, \(n_{\Hist}=8000\), \(n_\star=500\), \(\rho=0.5\), \(\tau=0.3\), and \(\theta_\star=0.5\), where \(n_{\Hist}\) counts all historical human and judge comparisons. Each setting uses \(500\) independent replications and ten cross-fitting folds. Supplementary Section~A.1 specifies the response-feature networks, coefficient scaling, and feature-shift construction.

\subsection{Estimation and inference}\label{subsec:sim-inference}

We report coverage of nominal \(95\%\) intervals, their mean width, and the root mean squared error (RMSE) of score estimates. ANCHOR uses an analytic variance estimator, avoiding repeated refitting for interval construction. Competing score estimators use basic bootstrap intervals from \(500\) full refits with shared resamples, accounting for historical estimation variability. HumanOnly uses a delta-method variance accounting for both human samples, and Oracle uses its arrival sandwich variance. Supplementary Section~A.8.3 gives the implementation details.

\begin{figure}[!htbp]
\centering
\IfFileExists{figures/simulation_interval_width_20260915/coverage_linear_main.pdf}{\includegraphics[width=\linewidth]{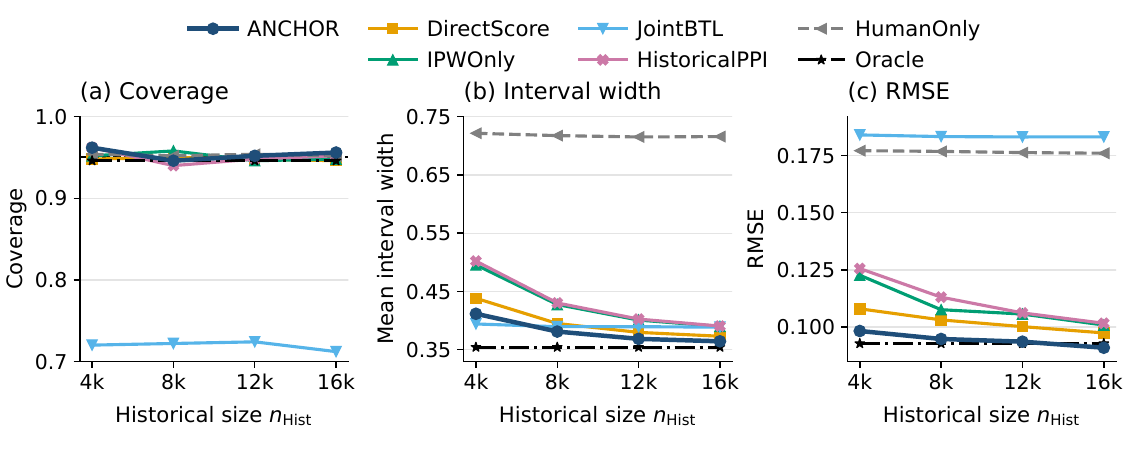}}{}
\caption{Coverage (a), mean interval width (b), and RMSE (c) as the historical sample grows, with \(n_\star=500\) and \(500\) replications per setting. The dashed line marks \(0.95\) coverage. Supplementary Figure~A.1 includes absolute bias and the full comparison of ten methods.}
\label{fig:sim-linear-validation}
\end{figure}

We vary \(n_{\Hist}\in\{4000,8000,12000,16000\}\) while holding \(n_\star=500\) fixed to examine the information gained from additional historical comparisons. ANCHOR's empirical coverage remains near nominal, and its RMSE approaches the Oracle value (Figure~\ref{fig:sim-linear-validation}). HistoricalPPI, DirectScore, and IPWOnly also have near-nominal coverage, with wider intervals and higher RMSE. At \(n_{\Hist}=8000\), RMSE is \(0.095\) for ANCHOR, \(0.113\) for HistoricalPPI, and \(0.177\) for HumanOnly. With \(150\) human comparisons involving the new model, HumanOnly has negligible empirical bias but higher variance than ANCHOR. JointBTL shows persistent bias and undercoverage, consistent with fitting a common comparison mechanism to heterogeneous human and judge outcomes. Supplementary Section~A.2 examines whether accounting for judge sensitivities and feature-dependent bias, averaging judge-specific estimates, or reweighting historical comparisons improves estimation and inference. Section~A.3 uses the same historical estimates to separately assess how weighting by judge sensitivity and the residual correction affect score estimation.

\begin{figure}[!htbp]
\centering
\IfFileExists{figures/variance_symbols_20260916/sparse_human_variance.pdf}{\includegraphics[width=\linewidth]{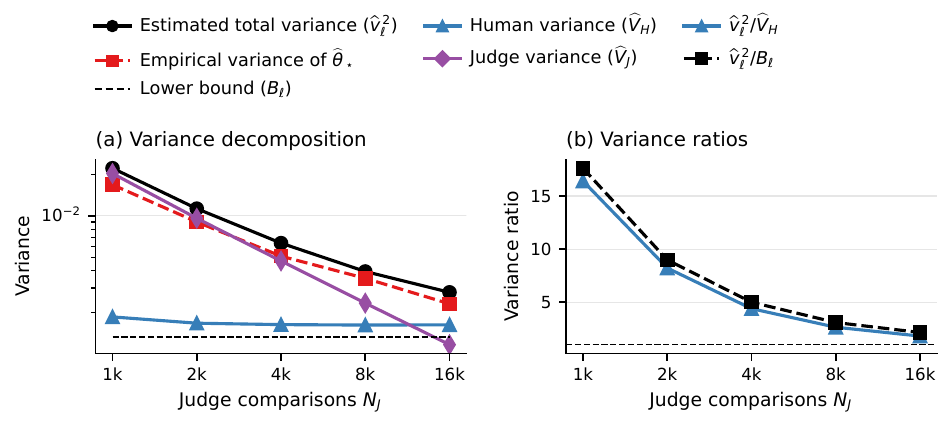}}{}
\caption{{Variance decomposition (a) and variance ratios (b), with \(n_{\Hum}=1000\) and \(\vect\ell=\vect e_1\). Variance estimates, components, and lower bounds are averaged over \(500\) replications. Ratios are computed within each replication and then averaged. Empirical variance is the sample variance of the \(500\) estimates of \(\theta_\star\).}}
\label{fig:sparse-human-variance}
\end{figure}

\label{subsec:sim-human-scale}
{
\clubpenalty=10000
We next examine how additional LLM-judge comparisons improve precision when the human sample size is fixed. For the new-model score \(T_\ell=\theta_\star\), with \(\vect\ell=\vect e_1\), we fix \(n_{\Hum}=1000\) and vary \(N_J\in\{1000,2000,4000,8000,16000\}\), allocating 70\% of LLM-judge comparisons to the historical stage and 30\% to the arrival stage. Each setting uses \(500\) independent replications. Here, \(\widehat v_\ell^2\) estimates the total variance of the corrected score estimator \(\widehat\theta_\star\). Its human and judge contributions are \(\widehat V_H\) and \(\widehat V_J\), respectively, with \(\widehat v_\ell^2=\widehat V_H+\widehat V_J\). The judge contribution includes both historical and arrival comparisons. Empirical variance is the sample variance of the \(500\) score estimates, and \(B_\ell\) is the theoretical variance lower bound from the human score-scale direction.

In Figure~\ref{fig:sparse-human-variance}(a), the average \(\widehat v_\ell^2\) is close to the empirical variance, and both decrease as \(N_J\) increases. The judge contribution \(\widehat V_J\) declines substantially, while the human contribution \(\widehat V_H\) changes little. The lower bound \(B_\ell\) depends only on the true scores and the fixed human comparison design, so its average does not change with \(N_J\). Figure~\ref{fig:sparse-human-variance}(b) shows that the average per-replication ratios \(\widehat v_\ell^2/\widehat V_H\) and \(\widehat v_\ell^2/B_\ell\) decline with \(N_J\) but remain above one at \(N_J=16000\). Thus, additional judge comparisons reduce estimation uncertainty, while the fixed human sample increasingly limits precision. Supplementary Section~A.1 gives the variance decomposition and ratio calculations.
}

\begin{figure}[!htbp]
\centering
\IfFileExists{figures/feature_shift_labels_20260915/two_step_bias_shift_grid.pdf}{\includegraphics[width=\linewidth]{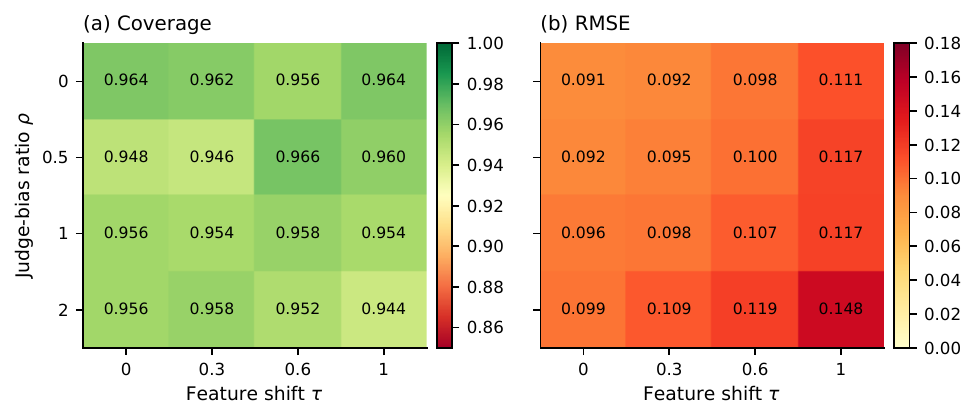}}{}
\caption{ANCHOR coverage (a) and RMSE (b) across judge-bias ratios \(\rho\) and feature shifts \(\tau\), with linear bias functions and the default sample sizes.}
\label{fig:bias-shift-grid}
\end{figure}

\label{subsec:sim-bias-shift}
We also examine the combined effects of judge bias and feature shift by varying \(\rho\in\{0,0.5,1,2\}\) and \(\tau\in\{0,0.3,0.6,1\}\) jointly. Coverage remains close to the nominal 95\% level across all settings, even as stronger judge bias and feature shift generally increase RMSE (Figure~\ref{fig:bias-shift-grid}). Supplementary Sections~A.4--A.6 examine sensitivity to the historical estimator, Riesz class, and comparison graph.

\subsection{Rank estimation}\label{subsec:sim-rank}

We examine how the gap between historical model scores affects the ranking accuracy of the new model. Historical scores are equally spaced and centered at zero, with adjacent gap \(\delta\in\{0.05,0.10,0.20,0.40\}\). The new model is evaluated at each of the nine adjacent-score midpoints. Judge sensitivities and bias coefficients remain fixed across gaps. All other settings follow the default design. Further details are given in Supplementary Section~A.1.

\begin{figure}[htb!]
\centering
\IfFileExists{figures/terminology_20260916/rank_insertion_gap.pdf}{\includegraphics[width=\linewidth]{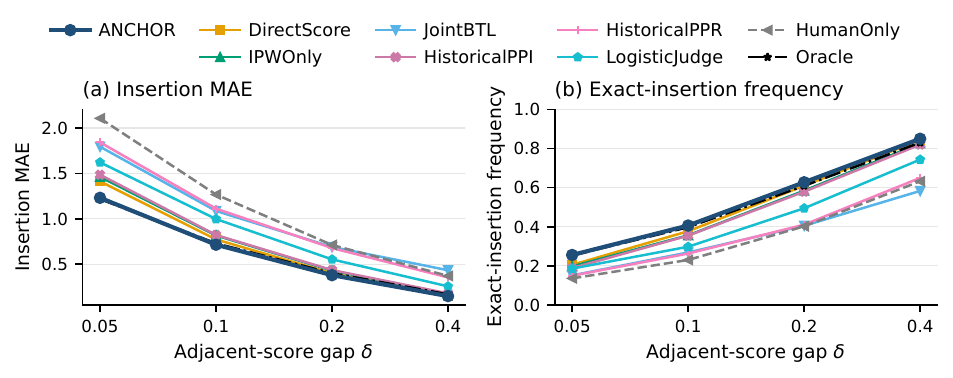}}{}
\caption{Insertion MAE (a) and exact-insertion frequency (b) for the seven comparison methods, HumanOnly, and Oracle as the gap \(\delta\) between adjacent historical scores increases.}
\label{fig:rank-insertion-gap}
\end{figure}

We fit a BTL model to the historical human comparisons to estimate the historical scores. Each method, including Oracle, compares its estimate for the new model with these historical scores to determine the insertion rank. For HistoricalPPR and LogisticJudge, we transform the historical scores using \textup{(A.3)} to make them comparable with the methods' ranking scores while preserving the historical ordering. Figure~\ref{fig:rank-insertion-gap}(a) reports insertion mean absolute error (MAE), the average absolute difference between estimated and true insertion ranks. Panel (b) reports exact-insertion frequency, the proportion of correctly estimated insertion ranks. Both metrics are averaged over the nine positions within each replication and then over the \(500\) replications.

{Apart from HumanOnly and Oracle, ANCHOR has the smallest observed insertion MAE at every gap in Figure~\ref{fig:rank-insertion-gap}(a). For example, at \(\delta=0.1\), its MAE is \(0.715\), about \(7.5\%\) lower than DirectScore's \(0.773\) and \(35.7\%\) lower than HistoricalPPR's \(1.112\). The exact-insertion frequencies are \(40.6\%\) for ANCHOR, \(37.7\%\) for DirectScore, and \(26.3\%\) for HistoricalPPR, giving ANCHOR gains of \(2.9\) and \(14.3\) percentage points, respectively. Thus, ANCHOR improves insertion accuracy relative to both the plug-in score estimator and the prediction-powered ranking comparator. As \(\delta\) increases from \(0.1\) to \(0.4\), ANCHOR's exact-insertion frequency rises from \(40.6\%\) to \(85.0\%\) (Figure~\ref{fig:rank-insertion-gap}(b)), consistent with insertion positions being easier to distinguish when adjacent historical scores are farther apart.} Supplementary Section~A.7 examines uncertainty in the new model's rank by reporting rank-set coverage and mean rank-set size together. It also considers model misspecification, under which high rank correlation can coexist with systematic insertion error.

\section{Application to Chatbot Arena}\label{sec:realdata}

The public Chatbot Arena data \citep{chiang2024chatbot} contain crowd preferences between model responses in user conversations. We obtain judgments from ten LLM judges for these response pairs. The judge request template is provided in Supplementary Section~B.2. The data cover 20 models and \(32{,}980\) prompts. After deduplication and removal of ties and unknown outcomes, \(211{,}209\) comparisons remain, including \(23{,}284\) human comparisons. We form a 12-dimensional feature vector for each response using principal components of its embedding. Supplementary Section~B.1 describes construction and preprocessing; Section~B.8 gives a second application to MT-Bench \citep{zheng2023judging}.

We design a leave-one-model-out (LOO) analysis in which each model is treated in turn as the new model, leaving \(K=19\) historical models. Historical fitting uses human and judge comparisons among these models; the held-out model's score is estimated using judge comparisons involving that model. Human comparisons involving the held-out model are reserved for evaluation. For insertion, each method compares its new-model estimate with historical scores from a human-only BTL fit to the 19 historical models. For evaluation, a separate human-only BTL fit to all 20 models provides the reference scores and ranks. We recenter these reference scores by subtracting their mean over the 19 historical models.

In each LOO analysis, we select seven judges based on their numbers of historical comparisons and distinct prompts, following the rule in Supplementary Section~B.3. All methods use the same selected judges. ANCHOR uses linear functions for judge bias and Riesz weight estimation, with ten-fold cross-fitting within each LOO analysis.

\begin{figure}[!htbp]
\centering
\includegraphics[width=\linewidth]{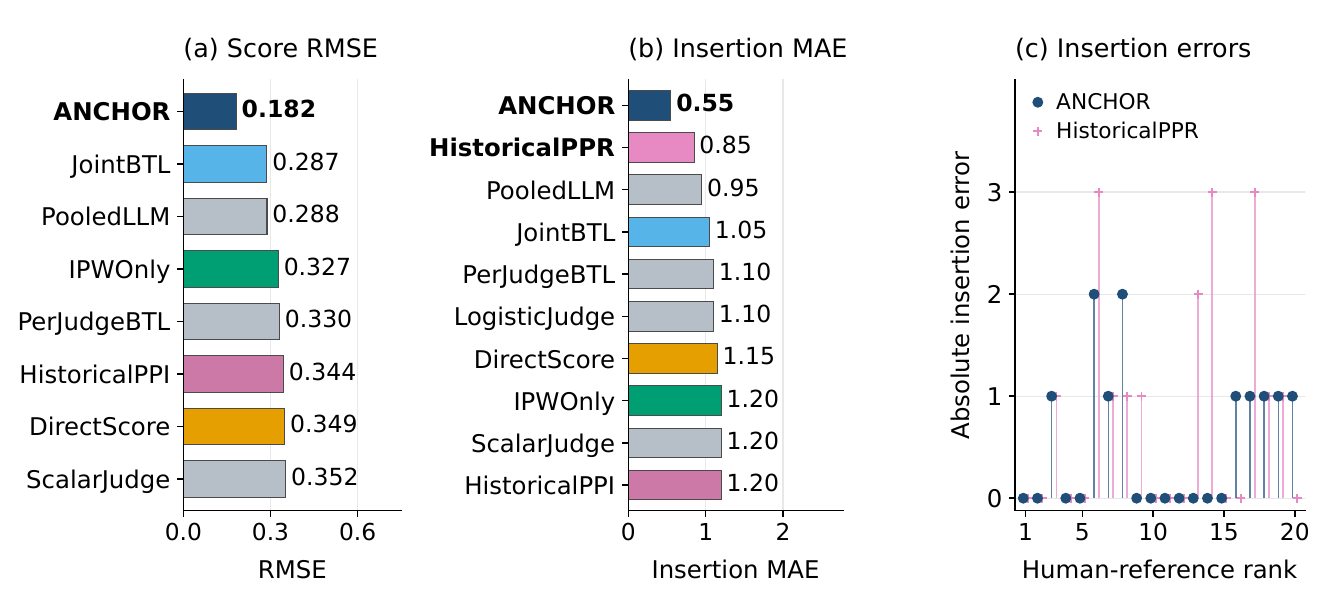}
\caption{Chatbot Arena score RMSE (a), insertion MAE (b), and absolute insertion errors for ANCHOR and HistoricalPPR (c), evaluated against the estimated human reference. Methods in panels (a) and (b) are ordered by error. In panel (c), the x-axis shows each model's human-reference rank, and the y-axis shows its absolute insertion error.}
\label{fig:real-loo}
\end{figure}

We assess agreement with the human benchmark through score RMSE and insertion errors. HistoricalPPR and LogisticJudge enter only the ranking comparison because their outputs are not on the human BTL score scale (Supplementary Section~A.8.3). Figure~\ref{fig:real-loo} summarizes score and insertion accuracy for the full set of competing methods evaluated in the main and supplementary simulations. Panel (c) compares ANCHOR with HistoricalPPR, which has the lowest insertion MAE among the other methods in panel (b). Supplementary Section~B.5 reports the complete comparisons, including Spearman correlations and exact-insertion counts in Table~B.3.

ANCHOR has the smallest score RMSE among the competing methods, \(0.182\), compared with \(0.287\) for JointBTL. Its insertion MAE is \(0.55\), approximately \(35\%\) lower than HistoricalPPR's \(0.85\). Supplementary Section~B.6 separates the contributions of likelihood weighting and orthogonal correction. The sensitivity analysis in Supplementary Section~B.7 shows that adding interactions between prompt and response features yields similar score estimates and interval widths and leaves all twenty predicted insertion ranks unchanged.

\begin{figure}[!htbp]
\centering
\includegraphics[width=\linewidth]{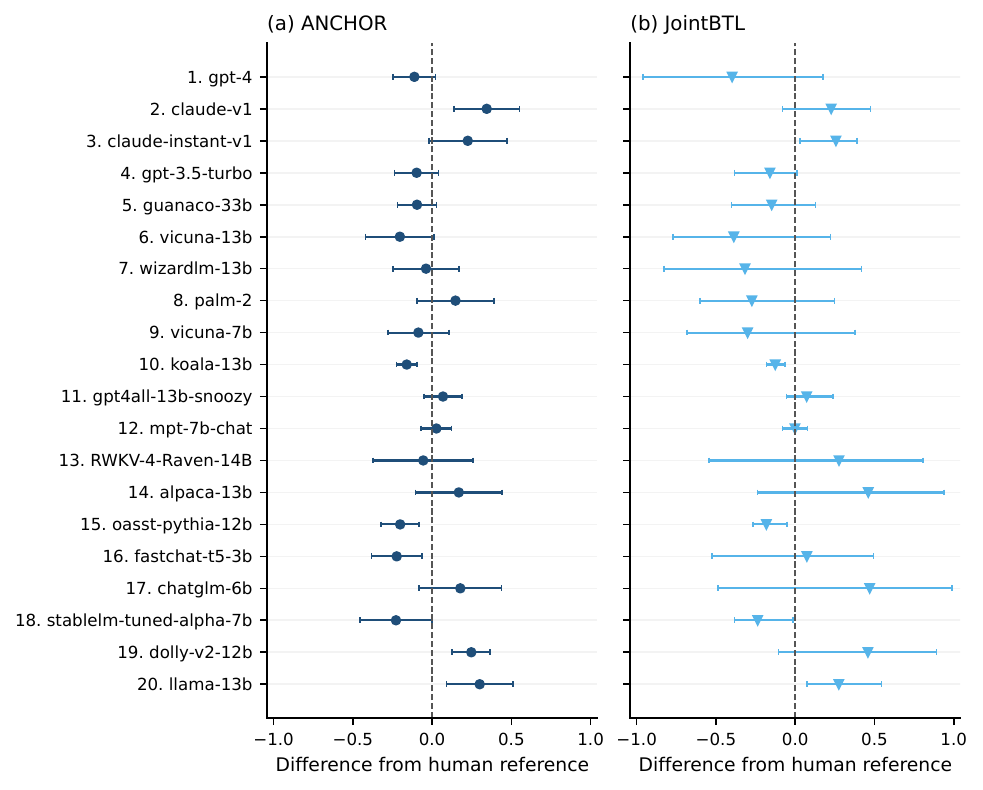}
\caption{Score estimates and nominal \(95\%\) intervals across 20 held-out Chatbot Arena models: analytic intervals for ANCHOR (a) and full-refit bootstrap intervals for JointBTL (b). Models are ordered by estimated human-reference rank.}
\label{fig:real-model-intervals}
\end{figure}

Figure~\ref{fig:real-model-intervals} compares ANCHOR's analytic intervals with the full-refit bootstrap intervals of JointBTL, the competing method with the lowest score RMSE. Supplementary Section~B.4 gives the interval constructions. Estimates and interval endpoints are centered at each model's estimated human-reference score, so zero marks the reference. ANCHOR's point estimates are closer to the human reference than JointBTL's for 14 of the 20 models. Its mean interval width is \(0.367\), compared with \(0.738\) for JointBTL, a reduction of approximately \(50\%\). The intervals contain the estimated human-reference scores for 14 of the 20 models, compared with 15 for JointBTL. Supplementary Section~B.9 shows that retaining all ten judges gives ANCHOR a similar score RMSE but lower human-reference inclusion. These comparisons describe agreement with an estimated benchmark: reference uncertainty is excluded from the intervals, the reference shares historical data with the fitted methods, and the LOO training samples overlap. Supplementary Table~B.2 reports reference inclusion and interval widths for all methods that estimate human-reference scores and for alternative interval constructions. Supplementary Table~B.4 gives ANCHOR's model-level scores, insertion ranks, and intervals alongside the human references.

\FloatBarrier
\section{Discussion}
\label{sec:discussion-v2}

\begingroup

ANCHOR establishes identification and efficient inference for new-model scores on a human-reference scale. Its orthogonal correction accounts for historical estimation error, while the efficiency analysis distinguishes the contributions of human and judge comparisons. The supplement presents additional theoretical extensions and experimental results. We briefly summarize these findings below and discuss their implications and future directions.

The theoretical extensions accommodate broader comparison structures. Supplementary Section~C.5 gives identification conditions for disconnected historical human comparison graphs, potentially involving only one pair of models with distinct scores. Section~D establishes identification and pointwise inference for scores and contrasts of fixed groups of new models connected internally by human comparisons and linked to historical models only through judge comparisons. Inference for joint historical fits with disconnected human comparisons remains a direction for future work.

Turning to estimation, joint historical fitting improves precision with small historical samples (Supplementary Section~A.4). Neural-network Riesz fits yield similar RMSE, wider intervals, and improved coverage and variance calibration compared with linear fits, although coverage remains below nominal for both classes (Supplementary Section~A.5). These findings motivate studying historical estimators and regularizing extreme correction weights.

Beyond estimation, the efficiency bound also guides selection of human comparison pairs and new-model opponents, and budget allocation. Supplementary Section~A.6 shows that sampling opponents from both groups of a sparsely connected historical graph reduces RMSE and interval width relative to sampling from one group at the same budget. Repeated updates motivate rank sets with simultaneous coverage that account for data reuse and adaptive selection.

Robustness to model misspecification remains an important question. The current guarantees allow feature shifts under shared judge sensitivities and bias functions. Supplementary Section~A.7 illustrates how omitted response-length preferences can shift insertion ranks while largely preserving estimated score ordering across settings (Table~A.5). Sensitivity analysis could bound the effects of such unrecorded preferences on score contrasts and identify when a comparison reverses. Human comparisons involving the new model could tighten these bounds.

A further extension would allow human preferences to depend on features, supporting context-specific rankings or overall rankings under a specified distribution of prompts and responses. Developing identification and orthogonal inference for these targets would extend human-anchored evaluation across tasks and deployment settings.
\endgroup

\begingroup
\setlength{\bibsep}{0pt}

{\fontsize{12}{14.4}\selectfont
\bibliography{bibliography}
}

\endgroup

\clearpage
\end{document}